\magnification = \magstep 1

\catcode`\@=11
\font\tensmc=cmcsc10      
\def\smc{\tensmc}

\def\hcorrection#1{\advance\hoffset by #1 }
\def\vcorrection#1{\advance\voffset by #1 }
\def\wlog#1{}
\newif\iftitle@
\outer\def\title{\title@true\vglue 24\p@ plus 12\p@ minus 12\p@
   \bgroup\let\\=\cr\tabskip\centering
   \halign to \hsize\bgroup\tenbf\hfill\ignorespaces##\unskip\hfill\cr}
\def\endtitle{\cr\egroup\egroup\vglue 18\p@ plus 12\p@ minus 6\p@}
\outer\def\author{\iftitle@\vglue -18\p@ plus -12\p@ minus -6\p@\fi\vglue
    12\p@ plus 6\p@ minus 3\p@\bgroup\let\\=\cr\tabskip\centering
    \halign to \hsize\bgroup\smc\hfill\ignorespaces##\unskip\hfill\cr}
\def\endauthor{\cr\egroup\egroup\vglue 18\p@ plus 12\p@ minus 6\p@}
\outer\def\heading{\bigbreak\bgroup\let\\=\cr\tabskip\centering
    \halign to \hsize\bgroup\smc\hfill\ignorespaces##\unskip\hfill\cr}
\def\endheading{\cr\egroup\egroup\nobreak\medskip}

\outer\def\proclaim#1{\medbreak\noindent\smc\ignorespaces
    #1\unskip.\enspace\sl\ignorespaces}
\outer\def\endproclaim{\par\ifdim\lastskip<\medskipamount\removelastskip
  \penalty 55 \fi\medskip\rm}
\outer\def\demo#1{\par\ifdim\lastskip<\smallskipamount\removelastskip
    \smallskip\fi\noindent{\smc\ignorespaces#1\unskip:\enspace}\rm
      \ignorespaces}
\outer\def\enddemo{\par\smallskip}
\newcount\footmarkcount@
\footmarkcount@=1
\def\makefootnote@#1#2{\insert\footins{\interlinepenalty=100
  \splittopskip=\ht\strutbox \splitmaxdepth=\dp\strutbox
  \floatingpenalty=\@MM
  \leftskip=\z@\rightskip=\z@\spaceskip=\z@\xspaceskip=\z@
  \noindent{#1}\footstrut\rm\ignorespaces #2\strut}}
\def\footnote{\let\@sf=\empty\ifhmode\edef\@sf{\spacefactor
   =\the\spacefactor}\/\fi\futurelet\next\footnote@}
\def\footnote@{\ifx"\next\let\next\footnote@@\else
    \let\next\footnote@@@\fi\next}
\def\footnote@@"#1"#2{#1\@sf\relax\makefootnote@{#1}{#2}}
\def\footnote@@@#1{$^{\number\footmarkcount@}$\makefootnote@
   {$^{\number\footmarkcount@}$}{#1}\global\advance\footmarkcount@ by 1 }

\hyphenation{man-u-script man-u-scripts ap-pen-dix ap-pen-di-ces}
\hyphenation{data-base data-bases}
\ifx\amstexloaded@\relax\catcode`\@=13
  \endinput\else\let\amstexloaded@=\relax\fi
\newlinechar=`\^^J
\def\eat@#1{}
\def\Space@.{\futurelet\Space@\relax}
\Space@. %
\newhelp\athelp@
{Only certain combinations beginning with @ make sense to me.^^J
Perhaps you wanted \string\@\space for a printed @?^^J
I've ignored the character or group after @.}
\def\futureletnextat@{\futurelet\next\at@}
{\catcode`\@=\active
\lccode`\Z=`\@ \lowercase
{\gdef@{\expandafter\csname futureletnextatZ\endcsname}
\expandafter\gdef\csname atZ\endcsname
   {\ifcat\noexpand\next a\def\next{\csname atZZ\endcsname}\else
   \ifcat\noexpand\next0\def\next{\csname atZZ\endcsname}\else
    \def\next{\csname atZZZ\endcsname}\fi\fi\next}
\expandafter\gdef\csname atZZ\endcsname#1{\expandafter
   \ifx\csname #1Zat\endcsname\relax\def\next
     {\errhelp\expandafter=\csname athelpZ\endcsname
      \errmessage{Invalid use of \string@}}\else
       \def\next{\csname #1Zat\endcsname}\fi\next}
\expandafter\gdef\csname atZZZ\endcsname#1{\errhelp
    \expandafter=\csname athelpZ\endcsname
      \errmessage{Invalid use of \string@}}}}
\def\atdef@#1{\expandafter\def\csname #1@at\endcsname}
\newhelp\defahelp@{If you typed \string\define\space cs instead of
\string\define\string\cs\space^^J
I've substituted an inaccessible control sequence so that your^^J
definition will be completed without mixing me up too badly.^^J
If you typed \string\define{\string\cs} the inaccessible control sequence^^J
was defined to be \string\cs, and the rest of your^^J
definition appears as input.}
\newhelp\defbhelp@{I've ignored your definition, because it might^^J
conflict with other uses that are important to me.}
\def\define{\futurelet\next\define@}
\def\define@{\ifcat\noexpand\next\relax
  \def\next{\define@@}%
  \else\errhelp=\defahelp@
  \errmessage{\string\define\space must be followed by a control
     sequence}\def\next{\def\garbage@}\fi\next}
\def\undefined@{}
\def\preloaded@{}
\def\define@@#1{\ifx#1\relax\errhelp=\defbhelp@
   \errmessage{\string#1\space is already defined}\def\next{\def\garbage@}%
   \else\expandafter\ifx\csname\expandafter\eat@\string
         #1@\endcsname\undefined@\errhelp=\defbhelp@
   \errmessage{\string#1\space can't be defined}\def\next{\def\garbage@}%
   \else\expandafter\ifx\csname\expandafter\eat@\string#1\endcsname\relax
     \def\next{\def#1}\else\errhelp=\defbhelp@
     \errmessage{\string#1\space is already defined}\def\next{\def\garbage@}%
      \fi\fi\fi\next}
\def\famzero{\fam\z@}

\def\lim{\mathop{\famzero lim}}

\def\textfont@#1#2{\def#1{\relax\ifmmode
    \errmessage{Use \string#1\space only in text}\else#2\fi}}
\textfont@\rm\tenrm
\textfont@\it\tenit
\textfont@\sl\tensl
\textfont@\bf\tenbf
\textfont@\smc\tensmc
\let\ic@=\/
\def\/{\unskip\ic@}
\def\textfonti{\the\textfont1 }
\def\t#1#2{{\edef\next{\the\font}\textfonti\accent"7F \next#1#2}}
\let\B=\=
\let\D=\.
\def~{\unskip\nobreak\ \ignorespaces}
{\catcode`\@=\active
\gdef\@{\char'100 }}
\atdef@-{\leavevmode\futurelet\next\athyph@}
\def\athyph@{\ifx\next-\let\next=\athyph@@
  \else\let\next=\athyph@@@\fi\next}
\def\athyph@@@{\hbox{-}}
\def\athyph@@#1{\futurelet\next\athyph@@@@}
\def\athyph@@@@{\if\next-\def\next##1{\hbox{---}}\else
    \def\next{\hbox{--}}\fi\next}
\def\.{.\spacefactor=\@m}
\atdef@.{\null.}
\atdef@,{\null,}
\atdef@;{\null;}
\atdef@:{\null:}
\atdef@?{\null?}
\atdef@!{\null!}
\def\srdr@{\thinspace}
\def\drsr@{\kern.02778em}
\def\sldl@{\kern.02778em}
\def\dlsl@{\thinspace}
\atdef@"{\unskip\futurelet\next\atqq@}
\def\atqq@{\ifx\next\Space@\def\next. {\atqq@@}\else
         \def\next.{\atqq@@}\fi\next.}
\def\atqq@@{\futurelet\next\atqq@@@}
\def\atqq@@@{\ifx\next`\def\next`{\atqql@}\else\def\next'{\atqqr@}\fi\next}
\def\atqql@{\futurelet\next\atqql@@}
\def\atqql@@{\ifx\next`\def\next`{\sldl@``}\else\def\next{\dlsl@`}\fi\next}
\def\atqqr@{\futurelet\next\atqqr@@}
\def\atqqr@@{\ifx\next'\def\next'{\srdr@''}\else\def\next{\drsr@'}\fi\next}

\def\textfontii{\the\textfont2 }
\def\{{\relax\ifmmode\lbrace\else
    {\textfontii f}\spacefactor=\@m\fi}
\def\}{\relax\ifmmode\rbrace\else
    \let\@sf=\empty\ifhmode\edef\@sf{\spacefactor=\the\spacefactor}\fi
      {\textfontii g}\@sf\relax\fi}
\def\nonhmodeerr@#1{\errmessage
     {\string#1\space allowed only within text}}
\def\linebreak{\relax\ifhmode\unskip\break\else
    \nonhmodeerr@\linebreak\fi}
\def\allowlinebreak{\relax
   \ifhmode\allowbreak\else\nonhmodeerr@\allowlinebreak\fi}
\newskip\saveskip@
\def\nolinebreak{\relax\ifhmode\saveskip@=\lastskip\unskip
  \nobreak\ifdim\saveskip@>\z@\hskip\saveskip@\fi
   \else\nonhmodeerr@\nolinebreak\fi}
\def\newline{\relax\ifhmode\null\hfil\break
    \else\nonhmodeerr@\newline\fi}
\def\nonmathaerr@#1{\errmessage
     {\string#1\space is not allowed in display math mode}}
\def\nonmathberr@#1{\errmessage{\string#1\space is allowed only in math mode}}
\def\mathbreak{\relax\ifmmode\ifinner\break\else
   \nonmathaerr@\mathbreak\fi\else\nonmathberr@\mathbreak\fi}
\def\nomathbreak{\relax\ifmmode\ifinner\nobreak\else
    \nonmathaerr@\nomathbreak\fi\else\nonmathberr@\nomathbreak\fi}
\def\allowmathbreak{\relax\ifmmode\ifinner\allowbreak\else
     \nonmathaerr@\allowmathbreak\fi\else\nonmathberr@\allowmathbreak\fi}
\def\pagebreak{\relax\ifmmode
   \ifinner\errmessage{\string\pagebreak\space
     not allowed in non-display math mode}\else\postdisplaypenalty-\@M\fi
   \else\ifvmode\penalty-\@M\else\edef\spacefactor@
       {\spacefactor=\the\spacefactor}\vadjust{\penalty-\@M}\spacefactor@
        \relax\fi\fi}
\def\nopagebreak{\relax\ifmmode
     \ifinner\errmessage{\string\nopagebreak\space
    not allowed in non-display math mode}\else\postdisplaypenalty\@M\fi
    \else\ifvmode\nobreak\else\edef\spacefactor@
        {\spacefactor=\the\spacefactor}\vadjust{\penalty\@M}\spacefactor@
         \relax\fi\fi}
\def\newpage{\relax\ifvmode\vfill\penalty-\@M\else\nonvmodeerr@\newpage\fi}
\def\nonvmodeerr@#1{\errmessage
    {\string#1\space is allowed only between paragraphs}}
\def\smallpagebreak{\relax\ifvmode\smallbreak
      \else\nonvmodeerr@\smallpagebreak\fi}
\def\medpagebreak{\relax\ifvmode\medbreak
       \else\nonvmodeerr@\medpagebreak\fi}
\def\bigpagebreak{\relax\ifvmode\bigbreak
      \else\nonvmodeerr@\bigpagebreak\fi}
\newdimen\captionwidth@
\captionwidth@=\hsize
\advance\captionwidth@ by -1.5in
\def\caption#1{}
\def\topspace#1{\gdef\thespace@{#1}\ifvmode\def\next
    {\futurelet\next\topspace@}\else\def\next{\nonvmodeerr@\topspace}\fi\next}
\def\topspace@{\ifx\next\Space@\def\next. {\futurelet\next\topspace@@}\else
     \def\next.{\futurelet\next\topspace@@}\fi\next.}
\def\topspace@@{\ifx\next\caption\let\next\topspace@@@\else
    \let\next\topspace@@@@\fi\next}
 \def\topspace@@@@{\topinsert\vbox to
       \thespace@{}\endinsert}
\def\topspace@@@\caption#1{\topinsert\vbox to
    \thespace@{}\nobreak
      \smallskip
    \setbox\z@=\hbox{\noindent\ignorespaces#1\unskip}%
   \ifdim\wd\z@>\captionwidth@
   \centerline{\vbox{\hsize=\captionwidth@\noindent\ignorespaces#1\unskip}}%
   \else\centerline{\box\z@}\fi\endinsert}
\def\midspace#1{\gdef\thespace@{#1}\ifvmode\def\next
    {\futurelet\next\midspace@}\else\def\next{\nonvmodeerr@\midspace}\fi\next}
\def\midspace@{\ifx\next\Space@\def\next. {\futurelet\next\midspace@@}\else
     \def\next.{\futurelet\next\midspace@@}\fi\next.}
\def\midspace@@{\ifx\next\caption\let\next\midspace@@@\else
    \let\next\midspace@@@@\fi\next}
 \def\midspace@@@@{\midinsert\vbox to
       \thespace@{}\endinsert}
\def\midspace@@@\caption#1{\midinsert\vbox to
    \thespace@{}\nobreak
      \smallskip
      \setbox\z@=\hbox{\noindent\ignorespaces#1\unskip}%
      \ifdim\wd\z@>\captionwidth@
    \centerline{\vbox{\hsize=\captionwidth@\noindent\ignorespaces#1\unskip}}%
    \else\centerline{\box\z@}\fi\endinsert}
\mathchardef\prime@="0230
\def\prime{{{}\prime@{}}}
\def\prim@s{\prime@\futurelet\next\pr@m@s}

\def\,{\relax\ifmmode\mskip\thinmuskip\else\thinspace\fi}
\def\!{\relax\ifmmode\mskip-\thinmuskip\else\negthinspace\fi}
\def\frac#1#2{{#1\over#2}}

\def\:{\nobreak\hskip.1111em{:}\hskip.3333em plus .0555em\relax}
\def\intic@{\mathchoice{\hskip5\p@}{\hskip4\p@}{\hskip4\p@}{\hskip4\p@}}
\def\negintic@
 {\mathchoice{\hskip-5\p@}{\hskip-4\p@}{\hskip-4\p@}{\hskip-4\p@}}
\def\intkern@{\mathchoice{\!\!\!}{\!\!}{\!\!}{\!\!}}
\def\intdots@{\mathchoice{\cdots}{{\cdotp}\mkern1.5mu
    {\cdotp}\mkern1.5mu{\cdotp}}{{\cdotp}\mkern1mu{\cdotp}\mkern1mu
      {\cdotp}}{{\cdotp}\mkern1mu{\cdotp}\mkern1mu{\cdotp}}}
\newcount\intno@
\def\iint{\intno@=\tw@\futurelet\next\ints@}
\def\iiint{\intno@=\thr@@\futurelet\next\ints@}
\def\iiiint{\intno@=4 \futurelet\next\ints@}
\def\idotsint{\intno@=\z@\futurelet\next\ints@}
\def\ints@{\findlimits@\ints@@}
\newif\iflimtoken@
\newif\iflimits@
\def\findlimits@{\limtoken@false\limits@false\ifx\next\limits
 \limtoken@true\limits@true\else\ifx\next\nolimits\limtoken@true\limits@false
    \fi\fi}
\def\multintlimits@{\intop\ifnum\intno@=\z@\intdots@
  \else\intkern@\fi
    \ifnum\intno@>\tw@\intop\intkern@\fi
     \ifnum\intno@>\thr@@\intop\intkern@\fi\intop}
\def\multint@{\int\ifnum\intno@=\z@\intdots@\else\intkern@\fi
   \ifnum\intno@>\tw@\int\intkern@\fi
    \ifnum\intno@>\thr@@\int\intkern@\fi\int}
\def\ints@@{\iflimtoken@\def\ints@@@{\iflimits@
   \negintic@\mathop{\intic@\multintlimits@}\limits\else
    \multint@\nolimits\fi\eat@}\else
     \def\ints@@@{\multint@\nolimits}\fi\ints@@@}
\def\Sb{_\bgroup\vspace@
        \baselineskip=\fontdimen10 \scriptfont\tw@
        \advance\baselineskip by \fontdimen12 \scriptfont\tw@
        \lineskip=\thr@@\fontdimen8 \scriptfont\thr@@
        \lineskiplimit=\thr@@\fontdimen8 \scriptfont\thr@@
        \Let@\vbox\bgroup\halign\bgroup \hfil$\scriptstyle
            {##}$\hfil\cr}
\def\endSb{\crcr\egroup\egroup\egroup}
\def\Sp{^\bgroup\vspace@
        \baselineskip=\fontdimen10 \scriptfont\tw@
        \advance\baselineskip by \fontdimen12 \scriptfont\tw@
        \lineskip=\thr@@\fontdimen8 \scriptfont\thr@@
        \lineskiplimit=\thr@@\fontdimen8 \scriptfont\thr@@
        \Let@\vbox\bgroup\halign\bgroup \hfil$\scriptstyle
            {##}$\hfil\cr}
\def\endSp{\crcr\egroup\egroup\egroup}
\def\Let@{\relax\iffalse{\fi\let\\=\cr\iffalse}\fi}
\def\vspace@{\def\vspace##1{\noalign{\vskip##1 }}}
\def\aligned{\,\vcenter\bgroup\vspace@\Let@\openup\jot\m@th\ialign
  \bgroup \strut\hfil$\displaystyle{##}$&$\displaystyle{{}##}$\hfil\crcr}
\def\endaligned{\crcr\egroup\egroup}
\def\matrix{\,\vcenter\bgroup\Let@\vspace@
    \normalbaselines
  \m@th\ialign\bgroup\hfil$##$\hfil&&\quad\hfil$##$\hfil\crcr
    \mathstrut\crcr\noalign{\kern-\baselineskip}}
\def\endmatrix{\crcr\mathstrut\crcr\noalign{\kern-\baselineskip}\egroup
                \egroup\,}
\newtoks\hashtoks@
\hashtoks@={#}
\def\format{\crcr\egroup\iffalse{\fi\ifnum`}=0 \fi\format@}
\def\format@#1\\{\def\preamble@{#1}%
  \def\c{\hfil$\the\hashtoks@$\hfil}%
  \def\r{\hfil$\the\hashtoks@$}%
  \def\l{$\the\hashtoks@$\hfil}%
  \setbox\z@=\hbox{\xdef\Preamble@{\preamble@}}\ifnum`{=0 \fi\iffalse}\fi
   \ialign\bgroup\span\Preamble@\crcr}
 
\let\hdots=\ldots
\def\cases{\left\{\,\vcenter\bgroup\vspace@
     \normalbaselines\openup\jot\m@th
       \Let@\ialign\bgroup$##$\hfil&\quad$##$\hfil\crcr
      \mathstrut\crcr\noalign{\kern-\baselineskip}}

\newif\iftagsleft@
\tagsleft@true
\def\TagsOnRight{\global\tagsleft@false}
\def\tag#1$${\iftagsleft@\leqno\else\eqno\fi
 \hbox{\def\pagebreak{\global\postdisplaypenalty-\@M}%
 \def\nopagebreak{\global\postdisplaypenalty\@M}\rm(#1\unskip)}%
  $$\postdisplaypenalty\z@\ignorespaces}
\interdisplaylinepenalty=\@M
\def\allowdisplaybreak@{\def\allowdisplaybreak{\noalign{\allowbreak}}}
\def\displaybreak@{\def\displaybreak{\noalign{\break}}}
\def\align#1\endalign{\def\tag{&}\vspace@\allowdisplaybreak@\displaybreak@
  \iftagsleft@\lalign@#1\endalign\else
   \ralign@#1\endalign\fi}
\def\ralign@#1\endalign{\displ@y\Let@\tabskip\centering\halign to\displaywidth
     {\hfil$\displaystyle{##}$\tabskip=\z@&$\displaystyle{{}##}$\hfil
       \tabskip=\centering&\llap{\hbox{(\rm##\unskip)}}\tabskip\z@\crcr
             #1\crcr}}
\def\lalign@
 #1\endalign{\displ@y\Let@\tabskip\centering\halign to \displaywidth
   {\hfil$\displaystyle{##}$\tabskip=\z@&$\displaystyle{{}##}$\hfil
   \tabskip=\centering&\kern-\displaywidth
        \rlap{\hbox{(\rm##\unskip)}}\tabskip=\displaywidth\crcr
               #1\crcr}}
\def\overrightarrow{\mathpalette\overrightarrow@}
\def\overrightarrow@#1#2{\vbox{\ialign{$##$\cr
    #1{-}\mkern-6mu\cleaders\hbox{$#1\mkern-2mu{-}\mkern-2mu$}\hfill
     \mkern-6mu{\to}\cr
     \noalign{\kern -1\p@\nointerlineskip}
     \hfil#1#2\hfil\cr}}}
\def\overleftarrow{\mathpalette\overleftarrow@}
\def\overleftarrow@#1#2{\vbox{\ialign{$##$\cr
     #1{\leftarrow}\mkern-6mu\cleaders\hbox{$#1\mkern-2mu{-}\mkern-2mu$}\hfill
      \mkern-6mu{-}\cr
     \noalign{\kern -1\p@\nointerlineskip}
     \hfil#1#2\hfil\cr}}}
\def\overleftrightarrow{\mathpalette\overleftrightarrow@}
\def\overleftrightarrow@#1#2{\vbox{\ialign{$##$\cr
     #1{\leftarrow}\mkern-6mu\cleaders\hbox{$#1\mkern-2mu{-}\mkern-2mu$}\hfill
       \mkern-6mu{\to}\cr
    \noalign{\kern -1\p@\nointerlineskip}
      \hfil#1#2\hfil\cr}}}
\def\underrightarrow{\mathpalette\underrightarrow@}
\def\underrightarrow@#1#2{\vtop{\ialign{$##$\cr
    \hfil#1#2\hfil\cr
     \noalign{\kern -1\p@\nointerlineskip}
    #1{-}\mkern-6mu\cleaders\hbox{$#1\mkern-2mu{-}\mkern-2mu$}\hfill
     \mkern-6mu{\to}\cr}}}
\def\underleftarrow{\mathpalette\underleftarrow@}
\def\underleftarrow@#1#2{\vtop{\ialign{$##$\cr
     \hfil#1#2\hfil\cr
     \noalign{\kern -1\p@\nointerlineskip}
     #1{\leftarrow}\mkern-6mu\cleaders\hbox{$#1\mkern-2mu{-}\mkern-2mu$}\hfill
      \mkern-6mu{-}\cr}}}
\def\underleftrightarrow{\mathpalette\underleftrightarrow@}
\def\underleftrightarrow@#1#2{\vtop{\ialign{$##$\cr
      \hfil#1#2\hfil\cr
    \noalign{\kern -1\p@\nointerlineskip}
     #1{\leftarrow}\mkern-6mu\cleaders\hbox{$#1\mkern-2mu{-}\mkern-2mu$}\hfill
       \mkern-6mu{\to}\cr}}}
\def\sqrt#1{\radical"270370 {#1}}
\def\dots{\relax\ifmmode\let\next=\ldots\else\let\next=\tdots@\fi\next}
\def\tdots@{\unskip\ \tdots@@}
\def\tdots@@{\futurelet\next\tdots@@@}
\def\tdots@@@{$\mathinner{\ldotp\ldotp\ldotp}\,
   \ifx\next,$\else
   \ifx\next.\,$\else
   \ifx\next;\,$\else
   \ifx\next:\,$\else
   \ifx\next?\,$\else
   \ifx\next!\,$\else
   $ \fi\fi\fi\fi\fi\fi}
\def\text{\relax\ifmmode\let\next=\text@\else\let\next=\text@@\fi\next}
\def\text@@#1{\hbox{#1}}
\def\text@#1{\mathchoice
 {\hbox{\everymath{\displaystyle}\def\textfonti{\the\textfont1 }%
    \def\textfontii{\the\textfont2 }\textdef@@ T#1}}
 {\hbox{\everymath{\textstyle}\def\textfonti{\the\textfont1 }%
    \def\textfontii{\the\textfont2 }\textdef@@ T#1}}
 {\hbox{\everymath{\scriptstyle}\def\textfonti{\the\scriptfont1 }%
   \def\textfontii{\the\scriptfont2 }\textdef@@ S\rm#1}}
 {\hbox{\everymath{\scriptscriptstyle}\def\textfonti{\the\scriptscriptfont1 }%
   \def\textfontii{\the\scriptscriptfont2 }\textdef@@ s\rm#1}}}
\def\textdef@@#1{\textdef@#1\rm \textdef@#1\bf
   \textdef@#1\sl \textdef@#1\it}

\def\textdef@#1#2{\def\next{\csname\expandafter\eat@\string#2fam\endcsname}%
\if S#1\edef#2{\the\scriptfont\next\relax}%
 \else\if s#1\edef#2{\the\scriptscriptfont\next\relax}%
 \else\edef#2{\the\textfont\next\relax}\fi\fi}
\scriptfont\itfam=\tenit \scriptscriptfont\itfam=\tenit
\scriptfont\slfam=\tensl \scriptscriptfont\slfam=\tensl
\mathcode`\0="0030
\mathcode`\1="0031
\mathcode`\2="0032
\mathcode`\3="0033
\mathcode`\4="0034
\mathcode`\5="0035
\mathcode`\6="0036
\mathcode`\7="0037
\mathcode`\8="0038
\mathcode`\9="0039
\def\Cal{\relax\ifmmode\let\next=\Cal@\else
     \def\next{\errmessage{Use \string\Cal\space only in math mode}}\fi\next}
\def\Cal@#1{{\fam2 #1}}
\def\bold{\relax\ifmmode\let\next=\bold@\else
   \def\next{\errmessage{Use \string\bold\space only in math
      mode}}\fi\next}\def\bold@#1{{\fam\bffam #1}}
\mathchardef\Gamma="0000
\mathchardef\Delta="0001
\mathchardef\Theta="0002
\mathchardef\Lambda="0003
\mathchardef\Xi="0004
\mathchardef\Pi="0005
\mathchardef\Sigma="0006
\mathchardef\Upsilon="0007
\mathchardef\Phi="0008
\mathchardef\Psi="0009
\mathchardef\Omega="000A
\mathchardef\varGamma="0100
\mathchardef\varDelta="0101
\mathchardef\varTheta="0102
\mathchardef\varLambda="0103
\mathchardef\varXi="0104
\mathchardef\varPi="0105
\mathchardef\varSigma="0106
\mathchardef\varUpsilon="0107
\mathchardef\varPhi="0108
\mathchardef\varPsi="0109
\mathchardef\varOmega="010A
\font\dummyft@=dummy
\fontdimen1 \dummyft@=\z@
\fontdimen2 \dummyft@=\z@
\fontdimen3 \dummyft@=\z@
\fontdimen4 \dummyft@=\z@
\fontdimen5 \dummyft@=\z@
\fontdimen6 \dummyft@=\z@
\fontdimen7 \dummyft@=\z@
\fontdimen8 \dummyft@=\z@
\fontdimen9 \dummyft@=\z@
\fontdimen10 \dummyft@=\z@
\fontdimen11 \dummyft@=\z@
\fontdimen12 \dummyft@=\z@
\fontdimen13 \dummyft@=\z@
\fontdimen14 \dummyft@=\z@
\fontdimen15 \dummyft@=\z@
\fontdimen16 \dummyft@=\z@
\fontdimen17 \dummyft@=\z@
\fontdimen18 \dummyft@=\z@
\fontdimen19 \dummyft@=\z@
\fontdimen20 \dummyft@=\z@
\fontdimen21 \dummyft@=\z@
\fontdimen22 \dummyft@=\z@
\def\fontlist@{\\{\tenrm}\\{\sevenrm}\\{\fiverm}\\{\teni}\\{\seveni}%
 \\{\fivei}\\{\tensy}\\{\sevensy}\\{\fivesy}\\{\tenex}\\{\tenbf}\\{\sevenbf}%
 \\{\fivebf}\\{\tensl}\\{\tenit}\\{\tensmc}}
\def\dodummy@{{\def\\##1{\global\let##1=\dummyft@}\fontlist@}}
\newif\ifsyntax@
\newcount\countxviii@
\def\newtoks@{\alloc@5\toks\toksdef\@cclvi}
\def\nopages@{\output={\setbox\z@=\box\@cclv \deadcycles=\z@}\newtoks@\output}
\def\syntax{\syntax@true\dodummy@\countxviii@=\count18
\loop \ifnum\countxviii@ > \z@ \textfont\countxviii@=\dummyft@
   \scriptfont\countxviii@=\dummyft@ \scriptscriptfont\countxviii@=\dummyft@
     \advance\countxviii@ by-\@ne\repeat
\dummyft@\tracinglostchars=\z@
  \nopages@\frenchspacing\hbadness=\@M}
\def\magstep#1{\ifcase#1 1000\or
 1200\or 1440\or 1728\or 2074\or 2488\or
 \errmessage{\string\magstep\space only works up to 5}\fi\relax}
{\lccode`\2=`\p \lccode`\3=`\t
 \lowercase{\gdef\tru@#123{#1truept}}}

\def\scaletype#1{\mag=#1\relax
 \hsize=\expandafter\tru@\the\hsize
 \vsize=\expandafter\tru@\the\vsize
 \dimen\footins=\expandafter\tru@\the\dimen\footins}

\def\scalefont#1#2\andcallit#3{\edef\font@{\the\font}#1\font#3=
  \fontname\font\space scaled #2\relax\font@}
\def\Mag@#1#2{\ifdim#1<1pt\multiply#1 #2\relax\divide#1 1000 \else
  \ifdim#1<10pt\divide#1 10 \multiply#1 #2\relax\divide#1 100\else
  \divide#1 100 \multiply#1 #2\relax\divide#1 10 \fi\fi}
\def\scalelinespacing#1{\Mag@\baselineskip{#1}\Mag@\lineskip{#1}%
  \Mag@\lineskiplimit{#1}}
\def\wlog#1{\immediate\write-1{#1}}
\catcode`\@=\active

\TagsOnRight

\title
A Generalization of Redfield's Master Theorem
\endtitle

\author
Valentin Vankov Iliev
\endauthor

\centerline{Section of Algebra, Institute of Mathematics,}
\centerline{Bulgarian Academy of Sciences, 1113 Sofia, Bulgaria}
\centerline{E-mail: viliev\@math.bas.bg,\ viliev\@nws.aubg.bg}

\heading
Introduction
\endheading

In our paper [3], we show that the natural environment for P\'olya's
fundamental enumeration theorem and for one of its possible generalizations, is
Schur-Macdonald's theory of invariant matrices (cf.  [4, Ch.  I, Appendix]).
The main result of this theory is the eqivalence between the category of
finite-dimensional linear representations of the symmetric group $S_d$ and the
category of polynomial homogeneous degree $d$ functors on the category of
finite-dimensional linear spaces.  The above-mentioned generalization is a
particular case of the main equality of characteristics of two objects
corresponding each other via Schur-Macdonald's equivalence.  For the original
version of P\'olya's theorem, we suppose that this fact was pointed out to G.
P\'olya by I. Schur himself (cf.  [5, Ch.  1, $n^o$ 20]).

Quite the contrary, it turns out that such a generalization of Redfield's
master theorem lives only in the environment of the representation theory of
the symmetric group, as a direct consequence of the decomposition of the tensor
product of several induced monomial representations of the symmetric group into
its transitive summands.  The underlying permutation representations give rise
to the original Redfield's group-reduced distributions, or, equivalently, to
Read's equivalence relation of \lq\lq $T$-similarity".

The paper is stratified as follows.  In Section 1, Proposition 1.2.5 asserts
that the tensor product of several induced monomial representations of the
symmetric group $S_d$ is a monomial one.  Lemma 1.2.3 discusses the
corresponding permutation representation of $S_d$ and, in particular, shows
that there is a canonical bijection between the $S_d$-orbit space thus obtained
and the set of Read's equivalence classes from [6, sec.  3].  Section 2
contains two equivalent statements:  Theorem 2.1.1 and Theorem 2.1.2, which
generalize Redfield's master theorem.  Evaluating the equality from Theorem
2.1.2 at $p_1=\cdots =p_d=1$, we get Theorem 2.1.3 which is a generalization of
the superposition theorem (cf.  [6, sec.  4]).  Using Lemma 1.2.3 as a link
between the representation theory of the symmetric group and the combinatorial
analysis, we formulate several particular cases of Theorem 2.1.3 as graphical
corollaries.  Let $\Gamma_1,\hdots ,\Gamma_k$ be graphs with $d$ vertices
(loops and multiple edges allowed), together with their automorphism groups
$W_1\leq S_d,\hdots , W_k\leq S_d$.  In Corollary 3.1.1 is given the number of
all superpositions of $\Gamma_1,\hdots ,\Gamma_k$, under the assumption that
the only $k$-tuple of permutations $(\sigma_1,\hdots ,\sigma_k)\in
W_1\times\cdots\times W_k$, each with the same cyclic structure, is
$((1),\hdots ,(1))$, where by $(1)$ we denote the unit permutation.  Suppose
that at least one of the groups $W_m$, say $W_1$, is cyclic and let $b$
be its order.  If $a$ is a divisor of $b$, the second corollary gives the
number of all superpositions of $\Gamma_1,\hdots ,\Gamma_k$ whose automorphism
groups are (cyclic) groups of order dividing $a$.  Under the same assumption,
Corollary 3.1.3 counts the number of all superpositions of $\Gamma_1,\hdots
,\Gamma_k$, whose automorphism groups are trivial.  Corollary 3.1.4 gives the
number of all superpositions of $\Gamma_1,\hdots ,\Gamma_k$, whose automorphism
groups consist of even permutations.

Below, with the kind permission of the participants of the e-mail discussion
apropos of the last two corollaries, we cite a fragment of it.

Prof. R. Robinson wrote me an e-letter where he, among other things, noted that
Corollary 3.1.3 \lq\lq\...  is a special case of Redfield's \lq\lq Case 2:
Theorem" , p.  449, where he counted the superpositions by automorphism group
when one of the range groups is cyclic". Prof. F. Harary concluded the
discussion with the following e-message:

\noindent \lq\lq Date: Fri, 8 Aug 1997 11:49:29 -0600 (MDT)

\noindent From: F. Harary $\langle\hbox{fnh\@crl.nmsu.edu}\rangle$

\noindent To: R. Robinson $\langle\hbox{rwr\@pollux.cs.uga.edu}\rangle$

\noindent Subject: V. Iliev's discovery

\noindent Cc: palmer\@math.msu.edu, rcread\@math.uwaterloo.ca,
viliev\@math.bas.bg

\noindent Dear Bob,

It is not at all surprising that Redfield anticipated Valentin as he also
anticipated Polya, Read, Ed, you, me and many many others!  I thank you very
much indeed for your prompt, detailed and scholarly reply to V. Iliev."

We note that [7, Case 2:  Theorem, p.  449] yields a finite procedure for
finding the number of superpositions with cyclic automorphism group of a given
order $a$ which divides the order of the cyclic range group, whereas in
Corollary 3.1.3 we give an explicit formula for the case $a=1$.

In case one of the automorphism groups is the dihedral group of order $2b$
where $b$ is an odd number, Corollary 3.2.1 establishes the number of all
superpositions of $\Gamma_1,\hdots ,\Gamma_k$, whose automorphism groups have
an odd order.

In general, when one of the automorphism groups has a normal solvable subgroup
of order $r$, such that the corresponding factor-group is cyclic of order
relatively prime to $r$, Corollary 3.2.2 counts the number of all
superpositions of $\Gamma_1,\hdots ,\Gamma_k$, whose automorphism groups have
an order dividing $r$.

Corollary 3.2.3 is a modification of 3.2.2.  Here the number $r$ is a power of
a prime $q$ and one of the $W_m$'s has a normal subgroup of order $r$, such
that the corresponding factor-group is cyclic of order relatively prime to $q$.
Under these hypotheses, the corollary gives the number of all superpositions of
$\Gamma_1,\hdots ,\Gamma_k$, whose automorphism groups are $q$-groups.

In Section 4 we lift the results from Theorem 2.1.1 and Theorem 2.1.2 via
Schur-Macdonald's equivalence to the category of polynomial homogeneous degree
$d$ functors on the category of finite-dimensional spaces.

\heading
1. Tensor product of induced monomial representations of $S_d$
\endheading

Throughout the end of the paper we assume that $K$ is an algebgaically closed
field of characteristic zero and that all group characters are $K$-valued.

1.1.  According to Schur-Macdonald's equivalence, the Grothendieck's group of
the category of polynomial homogeneous degree $d$\/ functors over the category
of finite dimensional $K$-linear spaces, is isomorphic to the Grothendieck's
group of the category of finite-dimensional $K$-linear representations of the
symmetric group $S_d$, that is, to the Abelian group $R^d$ consisting of all
generalized characters of $S_d$ (the set of all irreducible characters of $S_d$
is a basis for $R^d$).  The characteristic map $ch$ identifies both groups with
the Abelian group $\Lambda^d$ of homogeneous degree $d$\/ symmetric functions
with integer coefficients in a countable set of variables $x_0, x_1,
x_2,\hdots$ (cf.  [4, Ch.  I, Appendix, A7]).  The tensor product of two
finite-dimensional $K$-linear representations of $S_d$ with characters $u$ and
$v$, has character $uv$.  If $f=ch(u)$ and $g=ch(v)$, where $u$ and $v$ are
generalized characters of $S_d$, one defines the internal product $f\ast g$ of
two symmetric functions $f, g\in \Lambda^d$ by $f\ast g=ch(uv)$.  With respect
to the internal product the Abelian group $\Lambda^d$ becomes a commutative and
associative ring with identity element $h_d=ch(1_{S_d})$ [4, Ch.  I, sec.  7].
Transfering the product $*$ via the characteristic map $ch$ (or, eqiuvalently,
transfering the tensor product of representations of $S_d$ via
Schur-Macdonald's equivalence), one also defines {\it internal}
{\it product} $\ast$ of homogeneous degree $d$\/ polynomial
functors.

1.2. Let $W$ be a subgroup of the symmetric group $S_d$ and let $\chi\colon
W\to K$ be a one-dimensional character of $W$.
The field $K$ has a
natural structure of left $KW$-module given by $\sigma c=\chi (\sigma )c$,
where $\sigma\in W$, $c\in K$. We denote by $K_\chi$ the corresponding
one-dimensional $K$-linear representation of $W$.
Let $I$ be a left transversal
of $W$ in $S_d$.
The induced monomial representation
$ind_W^{S_d}(\chi )=KS_d\otimes_{KW}K_\chi$
has a natural basis $(e_i)_{i\in I}$, $e_i=i\otimes 1$, as a $K$-linear space.
Since for any
$\zeta\in S_d$ and $i\in I$ there exist unique $j\in I$ and $\sigma\in
W$ such that $\zeta i=j\sigma$, we obtain a group homomorphism
$s\colon S_d\to S(I)$
defined by the formula
$$
(s(\zeta)(i))^{-1}\zeta i\in W.\tag 1.2.1
$$
Moreover, the permutation group $s(S_d)$ is transitive on the set $I$. We have
$\zeta e_i=\zeta (i\otimes 1)=
(\zeta i)\otimes 1=(j\sigma)\otimes 1=j\otimes (\sigma 1)$.
Therefore
the
action of $S_d$ on
$ind_W^{S_d}(\chi )$ is given by
$$
\zeta e_i=
\beta_i(\zeta )
e_{s\left(\zeta\right)\left(i\right)},
$$
where
$\beta_i(\zeta )=\chi(\sigma)=\chi((s(\zeta)(i))^{-1}\zeta i)$.

For the rest of the paper we introduce the following notation.

$(W_m)_{m=1}^k$ is a finite family of subgroups
of the symmetric
group $S_d$;

$(\chi_m)_{m=1}^k$, $\chi_m\colon W_m\to K$, is a family
of one-dimensional characters;

For any
$m=1, 2,\hdots , k$, we denote by $I_m$,
$(e_i)_{i\in I_m}$,
$s_m\colon S_d\to S(I_m)$ and $(\beta_i^{\left(m\right)})_{i\in I_m}$, the
above ingredients for the induced monomial representation
$ind_{W_m}^{S_d}(\chi_m)$;

The rule
$(i_1,\hdots , i_k)\to (s_1(\zeta )(i_1),\hdots , s_k(\zeta )(i_k))$ defines a
group homomorphism
$$
s\colon S_d\to S(I),
i\to s(\zeta)(i), \tag 1.2.2
$$
where
$I=I_1\times\cdots\times I_k$ and
$i=(i_1,\hdots , i_k)$.

Next trivial lemma paves the way for some combinatorial applications.

\proclaim{Lemma 1.2.3} If
$s\colon S_d\to S(I)$ is the action (1.2.2)
of $S_d$ on the Cartesian product
$I=I_1\times\cdots\times I_k$, then

(i) two $k$-tuples
$(i_1,\hdots , i_k)$ and
$(j_1,\hdots , j_k)$ are in the same $S_d$-orbit in $I$ if and only if there
exist  $\zeta\in S_d$ and
$w_m\in W_m$, such that $j_m=\zeta i_mw_m$ for $m=1,\hdots , k$;

(ii) the stabilizer of the $k$-tuple
$(i_1,\hdots , i_k)$
in the symmetric group $S_d$, is the intersection
$i_1W_1i_1^{-1}\cap\hdots\cap i_kW_ki_k^{-1}$;

(iii) there exists a canonical bijection between the orbit space
$S_d\backslash I$ and the
orbit space
$$
S_d\times W_1^\circ\times\cdots\times W_k^\circ\backslash
S_d\times\cdots\times S_d,
$$
where the action of
the group $S_d\times W_1^\circ\times\cdots\times W_k^\circ$ on the set
$S_d\times\cdots\times S_d$ is given by
$$
(\zeta , w_1,\hdots , w_k)(a_1,\hdots , a_k)=(\zeta a_1w_1,\hdots
,\zeta a_kw_k)
$$
(here $W_m^\circ$ stands for the group
$W_m$ with the oposite group structure).

\endproclaim

\demo{Proof}  Definition of the action (1.2.2).

\proclaim{Remark 1.2.4}{\rm The
orbit space
$$
S_d\times W_1^\circ\times\cdots\times W_k^\circ\backslash
S_d\times\cdots\times S_d,
$$
coincides with the factor-set of
$S_d\times\cdots\times S_d$ with respect to the equivalence relation \lq\lq
$T$-similarity", defined in [6, sec. 3]. Therefore, according to Lemma 1.2.3,
(iii), there is a bijection between the orbit space $S_d\backslash I$ and the
set of
all distict superpositions of $k$ graphs with $d$ vertices each (multiple edges
and loops allowed), cf. [6, sec. 4].}

\endproclaim

\proclaim{Proposition 1.2.5}
The tensor product
$$
ind_{W_1}^{S_d}(\chi_1)\otimes_K
ind_{W_2}^{S_d}(\chi_1)\otimes_K\cdots\otimes_K
ind_{W_k}^{S_d}(\chi_k)\tag 1.2.6
$$
is a monomial $K$-linear representation of $S_d$ with basis
$(e_i=e_{i_1}\otimes\cdots\otimes e_{i_k})_{i\in I}$, the
action of $S_d$ being given by the rule
$$
\zeta e_i=\beta_i(\zeta )
e_{s\left(\zeta\right)(i)},
$$
where
$\beta_i(\zeta )=
\beta_{i_1}^{\left(1\right)}(\zeta )\hdots
\beta_{i_k}^{\left(k\right)}(\zeta ).$

\endproclaim

\demo{Proof}
It is clear that the family $(e_i)_{i\in I}$
is a basis for the $K$-linear space
(1.2.6). We have $\zeta e_i=\zeta e_{i_1}\otimes\cdots\otimes \zeta e_{i_k}=
\beta_{i_1}^{\left(1\right)}(\zeta )\hdots
\beta_{i_k}^{\left(k\right)}(\zeta )
e_{s_1\left(\zeta\right)(i_1)}\otimes\cdots\otimes
e_{s_k\left(\zeta\right)(i_k)}=
\beta_i(\zeta )
e_{s\left(\zeta\right)(i)}$.
In particular, (1.2.6) is a monomial
representation of $S_d$.

\enddemo

Let $N_0$ be the set of all non-negative integers. It is mentioned in [3, sec.
2] that the characteristic of an induced
monomial representation $ind_W^{S_d}(\chi )$ is equal to the generalized cyclic
index
$$
Z(\chi;p_1,\hdots , p_d)
=|W|^{-1}\sum_{\sigma\in W}\chi (\sigma)
p_1^{c_1\left(\sigma\right )}\hdots p_d^{c_d\left (\sigma\right )},
$$
where $p_s=\sum_{i\in N_0}x_i^s$ are the power sums and
$c_s(\sigma )$ is
the
number of cycles of length $s$ in the cyclic decomposition of the permutation
$\sigma$. By definition, the characteristic of the tensor product (1.2.6)
is the internal product
$$
Z(\chi_1;p_1,\hdots , p_d)\ast\cdots\ast Z(\chi_k;p_1,\hdots , p_d).
$$

We set
$$
C(W_1,\hdots , W_k)
=\{(\sigma_1,\hdots ,\sigma_k)\in W_1\times\cdots\times
W_k\mid c_s(\sigma_1)=\cdots =c_s(\sigma_k)\hbox{\ for\ }s=1,\hdots , d\}.
$$
Obviously, $((1),\hdots ,(1))\in C(W_1,\hdots , W_k)$.
For
$\sigma=(\sigma_1,\hdots ,\sigma_k)\in C(W_1,\hdots , W_k)$
we have
$\sigma^{-1}\in C(W_1,\hdots , W_k)$ and
$\eta\sigma\eta^{-1}\in C(W_1,\hdots , W_k)$ for each $\eta\in
W_1\times\cdots\times W_k$. In particular, the set
$C(W_1,\hdots , W_k)$ is a disjoint union of conjugasy classes of the abstract
group $W_1\times\cdots\times W_k$.

For any
$\sigma=(\sigma_1,\hdots ,\sigma_k)\in C(W_1,\hdots , W_k)$
we can define
$c_s(\sigma )=c_s(\sigma_1)=\cdots =c_s(\sigma_k)$ for $s=1,\hdots , d$.
Moreover, we set
$z_\sigma=\prod_{s=1}^ds^{c_s\left(\sigma\right)}c_s(\sigma)!$.

Next lemma links the present definition of internal product to Read's one from
[6, subsec. 3.3].

\proclaim{Lemma 1.2.7}
One has
$$
Z(\chi_1;p_1,\hdots , p_d)\ast\cdots\ast Z(\chi_k;p_1,\hdots , p_d)=
$$
$$
\frac{1}{|W_1|\hdots |W_k|}
\sum_{\sigma\in C\left(W_1,\hdots , W_k\right)}
z_\sigma^{k-1}\chi_1(\sigma_1)\hdots\chi_k(\sigma_k)
p_1^{c_1\left(\sigma\right )}\hdots p_d^{c_d\left (\sigma\right )}.
$$

\endproclaim

\demo{Proof} An immediate consequence of [4, Ch. I, sec. 7, (7.12)].

\enddemo

\heading
2. Redfield's Ansatz
\endheading

2.1. In this section we generalize Redfield's master theorem and
superposition theorem.

\proclaim{Theorem 2.1.1} One has
$$
ind_{W_1}^{S_d}(\chi_1)\otimes_K
ind_{W_2}^{S_d}(\chi_1)\otimes_K\cdots\otimes_K
ind_{W_k}^{S_d}(\chi_k)\simeq
$$
$$
\oplus_{\left(\omega_1,\hdots ,\omega_k\right)\in
T\left(W_1,\hdots , W_k\right)}
ind_{\omega_1W_1\omega_1^{-1}\cap\hdots\cap\omega_kW_k\omega_k^{-1}}^{S_d}
(\psi_{\left(\omega_1,\hdots ,\omega_k\right)}),
$$
where
$T\left(W_1,\hdots , W_k\right)$ is a system of distinct representatives of the
$S_d$-orbits in the Cartesian product $I=I_1\times\cdots\times I_k$ with
respect to the action (1.2.2) of $S_d$, and
$\psi_{\left(\omega_1,\hdots ,\omega_k\right)}$ is the one-dimensional
character of the group
$\omega_1W_1\omega_1^{-1}\cap\hdots\cap\omega_kW_k\omega_k^{-1}$, which is
the restriction of the expression
$$
\beta_\omega(\zeta)=\chi_1((s_1(\zeta)(\omega_1))^{-1}\zeta\omega_1)\hdots
\chi_k((s_k(\zeta)(\omega_k))^{-1}\zeta\omega_k)
$$
from (1.2.5).

\endproclaim

\demo{Proof} Proposition 1.2.5 implies that the tensor product (1.2.6) is a
monomial
representation of $S_d$, so it gives an induced monomial representation on
each $S_d$-orbit in the set $I$ and (1.2.6) is the direct sum of these
transitive
constituents. Now, Lemma 1.2.3, (ii), finishes the proof.

\enddemo

Transfering this result on the
Abelian group $\Lambda^d$ of homogeneous
degree $d$\/ symmetric functions with integer coefficients in a countable set
of variables  $x_0, x_1, x_2,\hdots$, we obtain a direct generalization of
Redfield's master theorem.

\proclaim{Theorem 2.1.2} One has
$$
Z(\chi_1;p_1,\hdots , p_d)\ast\cdots\ast Z(\chi_k;p_1,\hdots , p_d)=
$$
$$
\sum_{\left(\omega_1,\hdots ,\omega_k\right)\in
T\left(W_1,\hdots , W_k\right)}
Z(\psi_{\left(\omega_1,\hdots ,\omega_k\right)};p_1,\hdots , p_d).
$$

\endproclaim

Following R. C. Read, if $A$ is a polynomial in several variables $p_1,\hdots
, p_d$, we denote by $N(A)$ the sum of its coefficients.  Applying the operation
$N$ on the two sides of the previous equality, we establish

\proclaim{Theorem 2.1.3} One has that
$$
N(Z(\chi_1;p_1,\hdots , p_d)\ast\cdots\ast Z(\chi_k;p_1,\hdots , p_d))
$$
is the number of the elements $\omega\in
T\left(W_1,\hdots , W_k\right)$ such that
$\psi_{\left(\omega_1,\hdots ,\omega_k\right)}=1$ on the stabilizer
$\omega_1W_1\omega_1^{-1}\cap\hdots\cap\omega_kW_k\omega_k^{-1}$.

\endproclaim

\demo{Proof} It enough to note that given a
one-dimensional character
$\chi\colon W\to K$, we have
$N(Z(\chi;p_1,\hdots , p_d))=\langle\chi , 1_W\rangle_W$, where $\langle\hbox{\
,\ }
\rangle_W$ is the scalar product of functions on $W$, cf.
[4, Ch. I, sec. 7].

\enddemo

\proclaim{Remark 2.1.4}{\rm When $\chi_m=1_{W_m}$ for $m=1,\hdots , k$,
Theorem 2.1.2 (respectively, Theorem 2.1.3) turns into Redfield's master
theorem (respectively, turns into the superposition theorem).}

\endproclaim

\heading
Graphical corollaries
\endheading

3.1. Here is how the above
machinery applies to graph theory.
Combining Theorem 2.1.3 with Remark
1.2.4, we establish several graphical corollaries. Obviously, these corollaries
can
also be stated in the more general language of superpositions of $k$ sets each
of $d$ elements, with {\it a priori} given \lq\lq automorphism
groups".

Let $\Gamma_1,\hdots ,\Gamma_k$ be graphs with $d$
vertices (loops and multiple edges allowed) and let $W_1\leq
S_d,\hdots , W_k\leq S_d$, be their
automorphism groups, respectively.
Given a cyclic group of order $b$ and
a divisor $a$ of $b$,
let $\varrho^{\left(a\right)}$ be a one-dimensional character of this
cyclic group, whose  kernel has order $a$.

Lemma 1.2.7 yields immediately

\proclaim{Corollary 3.1.1} If $C(W_1,\hdots , W_k)=\{((1),\hdots ,(1))\}$, then
the number of all superpositions of
$\Gamma_1,\hdots ,\Gamma_k$ is
$$
\frac{(d!)^{k-1}}{|W_1|\hdots |W_k|}.
$$

\proclaim{Corollary 3.1.2} If the group $W_1$ is
cyclic of order $b$ and if $a$ is a divisor of $b$, then the number of all
superpositions
of $\Gamma_1,\hdots ,\Gamma_k$, which have a cyclic authomorphism group
of order dividing $a$, is
$$
N(Z(\varrho^{\left(a\right)};p_1,\hdots , p_d)\ast
Z(1_{W_2};p_1,\hdots , p_d)\ast
Z(1_{W_k};p_1,\hdots , p_d)),
$$
where $\varrho^{\left(a\right)}\colon W_1\to K$ ia a one-dimensional character
whose kernel has order $a$.

\endproclaim

In the particular case $a=1$, we obtain

\proclaim{Corollary 3.1.3} If the permutation group
$W_1$ is cyclic, then the number of
all superpositions of
$\Gamma_1,\hdots ,\Gamma_k$, having trivial authomorphism group, is
$$
N(Z(\varrho ;p_1,\hdots , p_d)\ast
Z(1_{W_2};p_1,\hdots , p_d)\ast
Z(1_{W_k};p_1,\hdots , p_d)),
$$
where $\varrho\colon W_1\to K$ is an injective one-dimensional character of
$W_1$.
\endproclaim

\proclaim{Corollary 3.1.4} The number of
all superpositions of
$\Gamma_1,\hdots ,\Gamma_k$, whose authomorphism group consist of even
permutations, is
$$
N(Z(\varepsilon;p_1,\hdots , p_d)\ast
Z(1_{W_2};p_1,\hdots , p_d)\ast
Z(1_{W_k};p_1,\hdots , p_d)),
$$
where $\varepsilon\colon W_1\to K$ is the restriction of the alternating
character of $S_d$ on $W_1$.

\endproclaim

3.2. Let us suppose that

(1) $W_1$ is the dihedral group of order $2b$ where $b$ is an odd
natural number.

Let $\delta\colon W_1\to K$ be the nontrivial one-dimensional
character of $W_1$. The kernel of $\delta$ is the cyclic subgroup $H\leq W_1$
of order
$b$ and all subgroups of $W_1$ which have an odd order, are subgroups of $H$.

\proclaim{Corollary 3.2.1}
If the group $W_1$ satisfies (1), then
the number of all superpositions
of $\Gamma_1,\hdots ,\Gamma_k$, which have an authomorphism group
of odd order, is
$$
N(Z(\delta ;p_1,\hdots , p_d)\ast
Z(1_{W_2};p_1,\hdots , p_d)\ast
Z(1_{W_k};p_1,\hdots , p_d)).
$$

\endproclaim

The last result is an important special case of the following statement.

Let $r$ be a natural number. We suppose that

(a) $W_1$ has a normal solvable subgroup $R$ of order $r$ such that the
factor-group $W_1/R$ is a cyclic group of order relatively prime to $r$.

Then the group $W_1$ itself is solvable.
According to the generalized Sylow theorems (cf [1, Ch. 9, Theorem
9.3.1]), $R$ is the only
subgroup of $W_1$ of order $r$. Moreover, any subgroup of $W_1$ of order
which divides $r$, is contained in $R$.

Denote by $\pi$ a one-dimensional character of $W_1$ with kernel $R$.

\proclaim{Corollary 3.2.2} If the group $W_1$ satisfies (a),
then the number
of all superpositions of $\Gamma_1,\hdots ,\Gamma_k$, whose authomorphism
groups are of order dividing $r$, coincide with
$$
N(Z(\pi ;p_1,\hdots , p_d)\ast
Z(1_{W_2};p_1,\hdots , p_d)\ast
Z(1_{W_k};p_1,\hdots , p_d)).
$$

\endproclaim

Now, we formulate a version of the preceding corollary. Let $q$ be a prime
number. Suppose that

(i) $W_1$ has a normal $q$-subgroup $R$ such that the factor-group $W_1/R$ is
a cyclic group of order relatively prime to $q$.

According to Sylow theorems (cf.  [1, Ch.  4, Theorems 4.2.1 - 4.2.3]), $R$ is
the only Sylow $q$-subgroup of $W_1$ and any $q$-subgroup of $W_1$ is contained
in $R$.

Denote by $\iota$ a one-dimensional character of $W_1$ with kernel $R$.

\proclaim{Corollary 3.2.3}
If the group $W_1$ satisfies (i), then
the number of all superpositions
of $\Gamma_1,\hdots ,\Gamma_k$, whose authomorphism
group is a $q$-group, coincide with
$$
N(Z(\iota ;p_1,\hdots , p_d)\ast
Z(1_{W_2};p_1,\hdots , p_d)\ast
Z(1_{W_k};p_1,\hdots , p_d)).
$$

\endproclaim

\proclaim{Remark 3.2.4}{\rm Examples of abstract groups $W_1$ which satisfy the
hypothesis (a) (respectively, the hypothesis (i)) can be obtained by
constructing
a semi-direct profuct of a solvable group $R$ of order $r$ (respectively, a
$q$-group $R$) with a cyclic group $C$ of order relatively prime to $r$
(respectively, relatively prime to $q$). Schur-Zassenhaus'
theorem
(cf.  [8, Ch.  IV, sec.  8, IV.7.c]) asserts that there are no other examples.
In the symmetric group $S_d$, it is enough to choose $R\leq
S_d$ and $C\leq
S_d$ with the above-mentioned properties so that $RC=CR$, $R\cap C=\{(1)\}$
and $R$ is a normal subgroup of the group $W_1=RC$.}

\endproclaim

\heading
4. Internal product of semi-symmetric powers
\endheading

4.1. Now, we transfer Theorem 2.1.1 as an isomorphism between polynomial
homogeneous degree $d$\/ functors, using Schur-Macdonald's equivalence. It is
known that the induced monomial representation
$ind_W^{S_d}(\chi )$
corresponds via that
equivalence to the polynomial functor $[\chi]^d(-)$ called a semi-symmetric
power (cf. [2, Theorem 2.2]).

\proclaim{Theorem 4.1.1} Using notation from (2.1.1), one has
$$
[\chi_1]^d(-)\ast
[\chi_1]^d(-)\ast\cdots\ast
[\chi_k]^d(-)\simeq
$$
$$
\oplus_{\left(\omega_1,\hdots ,\omega_k\right)\in
T\left(W_1\hdots , W_k\right)}
[\psi_{\left(\omega_1,\hdots ,\omega_k\right)}]^d(-).
$$

\endproclaim

The group
$\omega_1W_1\omega_1^{-1}\cap\hdots\cap\omega_kW_k\omega_k^{-1}\leq
S_d$ acts on the Cartesian product $N_0^d$ by permuting its components. Let
$J(N_0^d,\psi_{\left(\omega_1,\hdots ,\omega_k\right)})$
be the set of those
$\omega_1W_1\omega_1^{-1}\cap\hdots\cap\omega_kW_k\omega_k^{-1}$-orbits in
$N_0^d$, which contain the maximum number
$$
|\omega_1W_1\omega_1^{-1}\cap\hdots\cap\omega_kW_k\omega_k^{-1}:
Ker\psi_{\left(\omega_1,\hdots ,\omega_k\right)}|
$$
of $Ker\psi_{\left(\omega_1,\hdots ,\omega_k\right)}$-orbits (cf. [3, subsec.
3.2]). We set
$$
U(\chi_1,\hdots ,\chi_k)=
\coprod_{\left(\omega_1,\hdots ,\omega_k\right)\in
T\left(W_1\hdots , W_k\right)}
J(N_0^d,\psi_{\left(\omega_1,\hdots ,\omega_k\right)})
$$

\proclaim{Theorem 4.1.2}
The symmetric function
in a countable set of variables $x_0, x_1, x_2,\hdots$,
$$
Z(\chi_1;p_1,\hdots , p_d)\ast\cdots\ast Z(\chi_k;p_1,\hdots , p_d),
$$
is the weighted
inventory of the set
$U(\chi_1,\hdots ,\chi_k)$.

\endproclaim

\demo{Proof}
For any one-dimensional character $\chi\colon W\to K$
we denote by
$g(\chi,x_0,x_1,x_2,\hdots)$ the characteristics
$ch([\chi]^d(-)$ of the semi-symmetric power $[\chi]^d(-)$.
In view of Theorem 2.1.2 and [3, Theorem 2.1.2], we obtain that
$$
Z(\chi_1;p_1,\hdots , p_d)\ast\cdots\ast Z(\chi_k;p_1,\hdots , p_d)=
\sum_{\left(\omega_1,\hdots ,\omega_k\right)\in
T\left(W_1\hdots , W_k\right)}
g(\psi_{\left(\omega_1,\hdots ,\omega_k\right)};x_0,x_1,x_2,\hdots),
$$
and the statement follows from the combinatorial interpretation
of a symmetric function of the type
$g(\chi ;x_0,x_1,x_2,\hdots)$,
given in [3, subsec. 3.2].

\enddemo

\heading
References
\endheading

\noindent [1] M. Hall, Jr., The theory of groups, The Macmillan
Company, New York, 1959.

\noindent [2] V. V. Iliev, A note on the polynomial functors corresponding to
the monomial representations of the symmetric group, J. Pure and
Appl. Algebra 87 (1993), 1 -- 4.

\noindent [3] V. V. Iliev, A generalization of P\'olya's enumeration theorem or
the secret life of certain index sets, math/9902075

\noindent [4] I. G. Macdonald, Symmetric functions and Hall
polynomials, Clarendon Press, Oxford, 1995.

\noindent [5] G. P\'olya, Kombinatorische Anzahlbestimmungen f\"ur Gruppen,
Graphen und chemische Verbindungen, Acta Math. 68 (1937), 145 --
254. English translation:
G. P\'olya and R. C. Read, Combinatorial Enumeration of Groups,
Graphs and Chemical Compounds, Springer-Verlag New York Inc., 1987.

\noindent [6] R. C. Read, The enumeration of locally restricted graphs (I),
J. London Math. Soc. 34 (1959), 417 -- 436.

\noindent [7] J. H. Redfield, The theory of group-reduced distributions,
Amer. J. Math. 49 (1927), 433 -- 455.

\noindent [8] E. Schenkman, Group theory, D. Van Nostrand Company,
Inc., Princeton, New Jersey, 1965.

\end